\documentclass{amsart}
\usepackage{latexsym}
\usepackage[all]{xy}

\usepackage{amssymb} % \mathbb
\usepackage{amsmath} % \mathbb
\usepackage{amsthm}

\def\Z{{\mathbb{Z}}}% \Z == \mathbb{Z}
\def\R{{\mathbb{R}}}% \R == \mathbb{R}
\def\K{{\mathbb{K}}}% \K == \mathbb{K}
\def\A{{\mathcal{A}}}% \A == \mathcal{A}
\def\B{{\mathcal{B}}}% \B == \mathcal{B}
\DeclareMathOperator{\rank}{rank}

\DeclareMathOperator{\codim}{codim}

\DeclareMathOperator{\Der}{Der}

\numberwithin{equation}{section}

\newcommand{\owari}{\hfill$\square$}

\newtheorem{theorem}{Theorem}[section]
\newtheorem{prop}[theorem]{Proposition}
\newtheorem{cor}[theorem]{Corollary}
\newtheorem{lemma}[theorem]{Lemma}
\newtheorem{define}[theorem]{Definition}
\newtheorem{rem}[theorem]{Remark}

\newtheorem{example}[theorem]{Example}
\newtheorem{conjecture}[theorem]{Conjecture}

\newcommand{\graphA}[6]{
\begin{picture}(34,34)(-2,-2)
\xgraphA{#1}{#2}{#3}{#4}{#5}{#6}
\end{picture}
}
\newcommand{\GraphA}[7][1234]{
\raisebox{-20\unitlength}{
\begin{picture}(46,50)(-8,-8)
\xgraphAlabel#1
\xgraphA{#2}{#3}{#4}{#5}{#6}{#7}
\xgraphAvertex%
\end{picture}
}
}
\newcommand{\xgraphAlabel}[4]{
\put(-2,-2){\makebox(0,0)[tr]{\scriptsize $#3$}}
\put(32,-2){\makebox(0,0)[tl]{\scriptsize $#4$}}
\put(-2,32){\makebox(0,0)[br]{\scriptsize $#2$}}
\put(32,32){\makebox(0,0)[bl]{\scriptsize $#1$}}
}
\newcommand{\xgraphAvertex}[1][****]{
\xgraphAVertex #1
}
\newcommand{\xgraphAVertex}[4]{
\if#3o\put(0,0){\circle{4}}\fi%
\if#3*\put(0,0){\circle*{4}}\fi%
\if#3.\put(0,0){\circle*{2.4}}\fi%
\if#4o\put(30,0){\circle{4}}\fi%
\if#4*\put(30,0){\circle*{4}}\fi%
\if#4.\put(30,0){\circle*{2.4}}\fi%
\if#2o\put(0,30){\circle{4}}\fi%
\if#2*\put(0,30){\circle*{4}}\fi%
\if#2.\put(0,30){\circle*{2.4}}\fi%
\if#1o\put(30,30){\circle{4}}\fi%
\if#1*\put(30,30){\circle*{4}}\fi%
\if#1.\put(30,30){\circle*{2.4}}\fi%
}
\newcommand{\xgraphA}[6]{
\if#1+\put(30,30){\line(-1,0){30}}\fi% ~
\if#1.\qbezier[7](30,30)(15,30)(0,30)\fi% ~
\if#1-\put(30,31){\line(-1,0){30}}\put(30,29){\line(-1,0){30}}\fi% ~%%
\if#2+\put(0,0){\line(1,1){30}}\fi % /
\if#2.\qbezier[10](0,0)(15,15)(30,30)\fi % /
\if#2-\put(-0.7,0.7){\line(1,1){30}}\put(0.7,-0.7){\line(1,1){30}}\fi % /%%
\if#3+\put(30,30){\line(0,-1){30}}\fi% )
\if#3.\qbezier[7](30,0)(30,15)(30,30)\fi % /
\if#3-\put(29,30){\line(0,-1){30}}\put(31,30){\line(0,-1){30}}\fi% )%%
\if#4+\put(0,0){\line(0,1){30}}\fi % (
\if#4.\qbezier[7](0,0)(0,15)(0,30)\fi % /
\if#4-\put(-1,0){\line(0,1){30}}\put(1,0){\line(0,1){30}}\fi % (
\if#5+\put(0,30){\line(1,-1){30}}\fi % \
\if#5.\qbezier[10](30,0)(15,15)(0,30)\fi % /
\if#5-\put(-0.7,29.3){\line(1,-1){30}}\put(0.7,30.7){\line(1,-1){30}}\fi % \%%
\if#6+\put(0,0){\line(1,0){30}}\fi %_
\if#6.\qbezier[7](0,0)(15,0)(30,0)\fi % /
\if#6-\put(0,1){\line(1,0){30}}\put(0,-1){\line(1,0){30}}\fi %_
\xgraphAvertex}
\begin{document}

\title[Signed-eliminable graphs and free multiplicities]
{Signed-eliminable graphs and free multiplicities on the braid arrangement}
\author{Takuro Abe}
\address{Institute of Mathematics for Industry, Kyushu University, 
Motooka 744, Nishi-Ku, Fukuoka 8190395, Japan}
\email{abe@imi.kyushu-u.ac.jp}
\author{Koji Nuida}
\address{Information Technology Research Institute (ITRI), National Institute of Advanced Industrial Science and Technology (AIST), 2-3-26 Aomi, Koto-ku, Tokyo 1350064 Japan 
}
\email{k.nuida@aist.go.jp}
\author{Yasuhide Numata}
\address{Department of Mathematical Science, 
Faculty of Science, Shinshu University, Matsumoto, Nagano 3908621, Japan}
\email{nu@math.shinshu-u.ac.jp}



%\pagestyle{plain}

%%%%%%    TEXT START    %%%%%%

\maketitle

\begin{abstract}
We define specific multiplicities on the braid arrangement 
by using signed graphs. To consider their freeness, 
we introduce the notion of signed-eliminable graphs 
as a generalization of Stanley's 
classification theory of free graphic arrangements by chordal graphs. 
This generalization gives us a complete classification of the free 
multiplicities defined above. 
%We give a complete characterization of a certain class of 
%free multiplicities on the braid 
%arrangement induced from signed graphs 
%by introducing signed-eliminable graphs, 
%which are generalizations of 
%graphic arrangements and chordal graph theory by 
%Stanley.
%. This is a generalization of 
%that of graphic arrangements by Stanley, hence we introduce the concept of 
%signed-eliminable graphs as generalized chordal graphs.
%and its characterization is given. 
As an application, 
we prove one direction of a conjecture of 
Athanasiadis on the characterization of 
the freeness of 
certain deformations of 
the braid arrangement in terms of 
directed graphs. 
%in \cite{Ath2}.
\end{abstract}
\setcounter{section}{-1}

\section{Introduction}
Let $V=V^\ell$ be an $\ell$-dimensional vector space over a field $\K$ of 
characteristic zero,  
$\{x_1,\ldots,x_\ell\}$ a basis for the dual vector space $V^*$ and 
$S:=\mbox{Sym}(V^*) \simeq \K[x_1,\ldots,x_\ell]$. 
Let $\Der_{\K}(S)$ denote 
the $S$-module of $\K$-linear derivations of $S$, i.e., 
$
\Der_{\K}(S)=\bigoplus_{i=1}^{\ell} S \cdot \partial_{x_i}.
$
A non-zero element $\theta=\sum_{i=1}^\ell f_i \partial_{x_i} 
\in \Der_{\K}(S)$ is 
\textit{homogeneous of degree $p$} if $f_i$ is zero or homogeneous of 
degree $p$ for each $i$. 

A \textit{hyperplane arrangement} $\A$ (or simply 
an \textit{arrangement}) is a finite collection of affine hyperplanes 
in $V$. If each hyperplane in $\A$ contains the origin, we say that $\A$ is 
\textit{central}. In this article we assume that all arrangements are central 
unless otherwise specified. 
%The \textit{Coxeter arrangement of type $A_\ell$}, or equivalently, 
%the \textit{braid arrangement}, 
%sometimes denoted by $\A_\ell$, is defined as  
%$\{x_i-x_j=0|1 \le i<j\le \ell+1\}$ in $V^{\ell+1}$. 
A \textit{multiplicity} $m$ on an arrangement $\A$ is a map 
$m:\A \rightarrow \Z_{\ge0}$ and a pair $(\A,m)$ is called a 
\textit{multiarrangement}. 
Let $|m|$ denote the sum of the multiplicities 
$\sum_{H \in \A}m(H)$. When $m \equiv 1,\ (\A,m)$ is the same as the hyperplane arrangement 
$\A$ and sometimes called a \textit{simple arrangement}. For each hyperplane $H \in \A$ 
fix a linear form $\alpha_H \in V^*$ such that $\ker(\alpha_H)=H$. 
%Put 
%$
%Q(\A,m):=\prod_{H \in \A} \alpha_H^{m(H)}.
%$
The first main object in this article is the \textit{logarithmic derivation module} 
$D(\A,m)$ of $(\A,m)$ defined by 
$$
D(\A,m):=\{\theta \in \Der_{\K}(S)|\theta(\alpha_H) \in S\cdot \alpha_H^{m(H)}\ 
(\mbox{for all }H \in \A)\}.
$$
A multiarrangement $(\A,m)$ is \textit{free} if 
$D(\A,m)$ is a free $S$-module of rank $\ell$. If $(\A,m)$ is free, then 
there exists a homogeneous free basis $\{\theta_1,\ldots, \theta_\ell\}$ for 
$D(\A,m)$. Then we define the \textit{exponents} of a free multiarrangement 
$(\A,m)$ by $\exp(\A,m):=(\deg(\theta_1),\ldots,\deg(\theta_\ell))$. 
The exponents are independent of a 
choice of a basis. When $m \equiv 1$, 
%$(\A,m)$ is just an 
%arrangement $\A$ and we call it a \textit{simple arrangement}. In this case, 
the logarithmic derivation module and exponents are denoted by 
$D(\A)$ and $\exp(\A)$. When we fix a simple arrangement $\A$, we say that a 
multiplicity $m$ on $\A$ is \textit{free} (resp.$\ $\textit{non-free}) if 
a multiarrangement $(\A,m)$ is free (resp.$\ $non-free). 

A fundamental object of study in hyperplane arrangements is 
the arrangement of all reflecting hyperplanes of a Coxeter group, 
called a \textit{Coxeter arrangement}. The study of the 
logarithmic derivation module for a Coxeter arrangement and its freeness were  
%One of the main objects of study 
%in hyperplane arrangements is Coxeter arrangements and its 
%freeness, 
initiated by K. Saito in \cite{Sa1}, developed 
in \cite{Sa2}, and promoted by Solomon-Terao in \cite{ST2}, 
Terao in \cite{T4} and many other authors. 
In particular, Yoshinaga 
proved in \cite{Y2} and \cite{Y3} 
that the freeness of an arrangement is closely related to the
canonical restricted multiarrangement defined by Ziegler in \cite{Z}. 
% and he proved a conjecture of 
%Edelman and Reiner. 
Hence the freeness of multiarrangements 
is now a very important subject of research. 
%However, the structure of 
%$D(\A,m)$ is very difficult and there were only few approaches to it.

Recently, some results were 
%to consider $D(\A,m)$ are 
developed in \cite{ATW1} and \cite{ATW2} to study 
$D(\A,m)$ for general multiarrangements. 
%by Terao, Wakefield, 
%Yoshinaga, Yuzvinsky, 
%and the first author
% and third authors 
Also, some results concerning 
free multiplicities on Coxeter arrangements have been found, e.g., see 
\cite{A3}, \cite{AY2} and \cite{Y5}.
%, \cite{AN1}, 
%\cite{AN2}, \cite{WY} and \cite{Y5}. 
%By using these results, classifications of 
%free and non-free multiplicities are appearing by Terao, Wakefield, Yoshinaga and 
%the first author in \cite{ATW1}, \cite{ATW2}, \cite{Y5} and \cite{A3}. 
In this article we 
generalize the study of free multiplicities on the braid 
arrangement. 

A \textit{braid arrangement} $\A_\ell$, or the \textit{Coxeter arrangement of type $A_\ell$} is defined as 
$
\{H_{ij}:=\{x_i-x_j=0\}|1 \le i,j\le\ell+1,i \neq j\}
$
in $V=V^{\ell+1}$. 
By using the primitive derivation introduced in \cite{Sa1}, 
free multiplicities on Coxeter arrangements are studied by 
Solomon-Terao \cite{ST2}, 
Terao \cite{T4}, Yoshinaga \cite{Y0}, and the first author and Yoshinaga \cite{AY2}. 
%These three results are based on the theory of the primitive derivation by K. Saito 
%in \cite{Sa1}. 
Combining these results, 
we have a characterization of the freeness of quasi-constant multiplicities $m$ on 
a Coxeter arrangement, i.e., multiplicities such that $\max_{H,H' \in \A}|m(H)-m(H')| \le 1$. 
However, it is known that 
if $\max_{H,H' \in \A}|m(H)-m(H')| = 2$ then the 
same method using the primitive derivation 
does not work. Also, 
to determine explicitly which multiplicity makes $(\A,m)$ 
free is a difficult problem. 
%as the case of quasi-constant multiplicities. 
Our aim is 
to consider these multiplicities on the braid arrangement and 
classify their freeness completely. 
%, because the braid arrangement is 
%the simplest and fundamental Coxeter arrangement, and also because we can apply the addition-deletion 
%theorems for multiarrangements in \cite{ATW2} fully to the braid arrangement. Also, 
%the reason why 
In fact, we consider every multiplicity $m$ such that $|2k-m(H_{ij})| \le 1$ 
for some $k \in \Z_{>0}$ 
since, as shown in \cite{AY2}, a mysterious and interesting symmetry 
of the freeness and duality of exponents exists for these kinds of multiplicities $m$. 
%Hence we can apply the original theory of free hyperplane arrangements and 
%graph theory to these kinds of multiplicities, 
%and cannot to those which has an odd constant multiplicity as its center, 
%see \cite{ATW1} and \cite{AY2}.

To state the main theorem, let us introduce some notation. Let $\A$ be the 
braid arrangement in $V^{\ell+1}$. To express the multiplicity $m$ mentioned in the 
previous paragraph, 
%satisfying  
%$\max_{H,H' \in \A}|m(H)-m(H')| \le 2$, 
we use a 
signed graph $G$, i.e., $G$ is a graph consisting of 
the vertex set $V_G=\{v_1,v_2,\ldots,v_{\ell+1}\}$ and 
the set of edges $E_G$ which has the decomposition $E_G=E_G^+ \cup E_G^-$ with 
$E_G^+ \cap E_G^-=\emptyset$. Then we can define the following 
map. 

\begin{define}\label{multiplicity}
The map $m_G$ on 
the braid arrangement $\A_\ell$ is defined by 
\[
m_G(H_{ij}):=
\left\{
\begin{array}{rl}
1 & \mbox{if}\ \{v_i,v_j\} \in E_G^+ ,\\
-1 & \mbox{if}\ \{v_i,v_j\} \in E_G^-,\ \mbox{and}\\
0 & \mbox{otherwise}, 
\end{array}
\right.
\]
where $\{v_i,v_j\}$ denotes the undirected edge between $v_i$ and $v_j$. 
\end{define}

Also, we introduce the following notion 
of signed graphs to characterize 
the freeness. 

%To characterize the freeness of $m_G$, let us introduce the following 
%concept 
%of
%graphs.

\begin{define}
The graph $G$ is \textit{signed-eliminable} with a \textit{signed-elimination ordering} 
$\nu:V_G \rightarrow \{1,2,\ldots,\ell+1\}$ if $\nu$ is bijective, and 
for every three vertices $v_i,v_j,v_k \in V_G$ with $\nu(v_i),\nu(v_j)<\nu(v_k)$, the induced subgraph 
$G|_{\{v_i,v_j,v_k\}}$ satisfies the following conditions:
\begin{itemize}
\item[(1)] For $\sigma \in \{+,-\}$, 
if $\{v_i,v_k\}$ and $\{v_j,v_k\}$ are edges in $E_G^\sigma$, then 
$\{v_i,v_j\} \in E_G^\sigma$.
%\item[(2)] $\{v_i,v_k\}$ and $\{v_j,v_k\}$ are edges of the same color, and $\{v_i,v_j\}$ is an edge of 
%the different color.
\item[(2)] For $\sigma \in \{+,-\}$, 
if $\{v_k,v_i\} \in E_G^\sigma$ and $\{v_i,v_j\} \in E_G^{-\sigma}$, then 
$\{v_k,v_j\} \in E_G$.
\end{itemize}
For a signed-eliminable graph $G$ with a signed-elimination ordering 
$\nu$, $v \in V_G$ and $i \in \{1,2,\ldots,\ell+1\}$, 
define the degree $\widetilde{\deg}_{i}(v)$ by 
$$
\widetilde{\deg}_{i}(v):=\deg(v,V_G, E_G^+|_{\nu^{-1}\{1,2,\ldots,i\}})
-\deg(v,V_G, E_G^-|_{\nu^{-1}\{1,2,\ldots,i\}}),
$$
where $\deg(w,V_H,E_H):=|\{x \in V_H|\{w,x\} \in E_H\}|$ is the degree of the vertex $w$ in 
the graph $H=(V_H,E_H)$, and  
$(V_G,E_G^{\sigma}|_S)$ with respect to $S \subset V_G$ is the induced subgraph of $G$ whose 
set of edges is equal to $\{\{v_i,v_j\} \in E_G^\sigma|v_i,v_j \in S\}$. 
Furthermore, define $\widetilde{\deg}_i:=
\widetilde{\deg}_i(\nu^{-1}(i))$ for each $i\ (1 \le i \le \ell+1)$. 
\label{bichordal}
\end{define}

We consider the property of signed-eliminable graphs in Sections two and three. 
Also note that a signed-eliminable graph is a generalization of a chordal graph, or 
a graph which has a vertex elimination order (see Remark \ref{remark}). 
By using chordal graphs, Stanley classified completely 
the free and non-free graphic arrangements in \cite{S} (see also \cite{ER1} or 
Section one in this article). 
What we will do in this article is the multi-version of Stanley's result. 
In other words,  
%His theory gives a very explicit and beautiful characterization of 
%free arrangements. 
%As we mentioned, there exists a shifting isomorphism which characterizes 
%the freeness of $(\A,2k+m)$ for $m:\A \rightarrow \{0,1\}$ and $k \in \Z_{>0}$ in 
%\cite{AY2}. 
%However, to determine explicitly 
%which multiarrangement $(\A,m)$ is free or not free, is a totally different problem. 
%From this viewpoint,  
we will classify free multiplicities 
on the braid arrangement of the form $2k+m_G$
%do the same classification for the multi-braid arrangement with 
%multiplicities $2k+m_G$ for $m_G$ 
with $m_G$ defined in Definition \ref{multiplicity} in more general setting.  
The main result is the following characterization of the 
freeness 
%multiplicities on the braid arrangement 
in terms of signed-eliminable graphs.\footnote{
Theorem \ref{freemulti} is not correct as it is stated below. 
See Appendix A for the corrected statements, conditions and proofs.}

\begin{theorem}
Let $\A$ be the braid arrangement in $V^{\ell+1}$, $G$ a signed graph 
and $m_G$ the map in Definition \ref{multiplicity}.
Let $k,
n_1,\ldots,n_{\ell+1}$ be non-negative integers. Define a multi-braid arrangement 
$(\A,m)=\A_\ell(n_1,n_2,\ldots,n_{\ell+1})[G]$ by 
%$\A :=\{H_{ij}=\{x_i-x_j=0\}|1 \le i<j\le \ell+1\}$ and 
$m(H_{ij})=2k+n_i+n_j+m_G(H_{ij})$ and 
put $N=(\ell+1)k+\sum_{i=1}^{\ell+1} n_i$. 
Assume that one of the following three 
conditions is satisfied: 
\begin{itemize}
\item[(a)] 
$k>0$.
\item[(b)] 
$E_G^- = \emptyset$.
\item[(c)] 
$E_G^+ = \emptyset$ and $m(H_{ij})>0$ for all $H_{ij} \in \A$. 
\end{itemize}
Then 
$\A_\ell(n_1,n_2,\ldots,n_{\ell+1})[G]$ is free with  
$$
\exp(\A,m)=(0,N+\widetilde{\deg}_{2},\ldots,N+\widetilde{\deg}_{{\ell+1}})
$$
if and only if 
$G$ is signed-eliminable.
% with a signed-elimination ordering $\nu$.
\label{freemulti}
\end{theorem}

If we let $n_1=\cdots=n_{\ell+1}=0$ for case (b) of Theorem 0.3 
then the corresponding arrangement is a graphic arrangement where each hyperplane has multiplicity one. 
Therefore, Theorem \ref{freemulti} is a generalization of Stanley's classification of free graphic arrangements. 
In Sections two and three we will see that a signed-eliminable graph is 
a generalization of the concept of a chordal graph. 
Hence, Theorem \ref{freemulti} 
generalizes both aspects of Stanley's work in \cite{S}: 
the freeness of certain arrangements and combinatorial properties of the corresponding graphs.

%The classification of free graphic arrangements corresponds to the 
%case (b) in Theorem \ref{freemulti} with $n_1=\ldots=n_{\ell+1}=0$. 
%As mentioned before, we will see later that a signed-eliminable graph is 
%a generalization of a chordal graph. Hence Theorem \ref{freemulti} 
%generalizes both aspects of Stanley's theory in \cite{S}, i.e., 
%the freeness of arrangements and combinatorics of graphs. 
%%Hence this article is a generalization of \cite{S}. 

The organization of this article is as follows. In Section one 
we introduce some fundamental results and definitions about multiarrangements and 
their freeness. 
In Section two we introduce the theory of signed-eliminable graphs, which can be 
regarded as a generalization of the chordal graph theory from the viewpoint of 
the characterization of free graphic arrangements due to Stanley. 
In Section three we quote a characterization of 
signed-eliminable graphs from \cite{N}. 
In Section four we apply the results in the previous sections to 
the study of free multiplicities on the braid arrangement, and prove 
Theorem \ref{freemulti}. 
In Section five, we give an application of Theorem \ref{freemulti} 
to a conjecture of Athanasiadis in \cite{Ath2}. 
\medskip

\textbf{Acknowledgments}. The authors appreciate Professor 
Masahiko Yoshinaga and Professor 
Max Wakefield for advice and comments to this article. 
Also 
the authors are grateful to the referee for useful comments to this article. 
The first and third authors were supported by 21st Century COE Program 
``Mathematics of Nonlinear Structures via Singularities'' Hokkaido University.
%on the 
%duality of exponents in Theorem \ref{freemulti}. 

\section{Preliminaries}
In this section let us review some results and definitions which 
will be used in this article. Let us begin with those for (multi)arrangements of hyperplanes, 
for which we refer the reader to \cite{OT}. 
First we introduce some results for the study of free and non-free multiarrangements. 
Let $(\A,m)$ be a multiarrangement in an $\ell$-dimensional vector space 
and fix $H_0 \in \A$ with $m(H_0)>0$. 
Define the \textit{deletion} $(\A',m')$ of $(\A,m)$ 
with respect to $H_0$ 
by $\A'=\A$ and 
\[
m'(H)=
\left\{
\begin{array}{rl}
m(H) \ \ \ \ \ \  & \mbox{if}\ H \neq H_0 ,\\
m(H_0)-1 & \mbox{if}\ H = H_0.
\end{array}
\right.
\]
%Hence for the deletion $(\A',m')$ 
%the arrangement $\A'$ is equal to $\A \setminus \{H_0\}$ or $\A$ depending on 
%$m(H_0)=1$ or not.

\begin{theorem}[\cite{ATW2}, Theorem 0.4]
If $(\A,m)$ and $(\A',m')$ are both free, then 
there exists a basis $\{\theta_1,\ldots,\theta_\ell\}$ for $D(\A',m')$ such that 
$\{\theta_1,\ldots,\theta_{k-1}, 
\alpha_{H_0}\theta_k,\theta_{k+1},\ldots,\theta_\ell\}$ is a basis for $D(\A,m)$ 
for some $k \in \{1,\ldots, \ell\}$.
\label{basis}
\end{theorem}

For $X \in \A'':= \{H' \cap H_0|H' \in \A \setminus \{H_0\}\}$, 
define $\A_X:=\{H \in \A|X \subset H\}$ and $m_X:=m|_{\A_X}$. 
Since $\A_X$ is essentially a $2$-multiarrangement, Theorem \ref{basis} implies that 
$(\A_X,m_X)$ is free with a basis 
$\{\zeta_3,\zeta_4,\ldots,\zeta_\ell,\theta_X,\psi_X\}$, where 
$\deg(\zeta_i)=0$, $\theta_X \not \in \alpha_{H_0} \mbox{Der}_{\K}(S) $ and 
$\psi_X \in \alpha_{H_0} \mbox{Der}_{\K}(S) $. 
%(both mod $\alpha_{H_0}$)
Then we define the 
\textit{Euler multiplicity} $m^*$ on $\A''$ by 
$m^*(X):=\deg(\theta_X)$, and we call $(\A'',m^*)$ the \textit{Euler restriction}. 
Then the following \textit{Addition-Deletion theorem} holds.

%\begin{theorem}[\cite{ATW2}, Theorem 0.4]
%If a multiarrangement $(\A,m)$ is free, then 
%$GMP_k=LMP_k\ (1 \le k \le \ell)$, where $GMP_k$ is the $k$-th global mixed product of 
%$(\A,m)$ and 
%$LMP_k$ is the $k$-th local mixed product of $(\A,m)$.
%\label{gmplmp}
%\end{theorem}

\begin{theorem}[\cite{ATW2}, Theorem 0.8]
%Let $(\A,m)$ be a nonempty multiarrangement in an $\ell$-dimensional 
%vector space $V$, $H_0 \in \A$ 
%and 
Let $(\A,m),(\A',m')$ and $(\A'',m^{*})$ be the triple with respect to $H_0$. Then any two of the following 
statements imply the third:
\begin{itemize}
\item[(i)]
$(\A,m)$ is free with $\exp(\A,m)=(d_1,\ldots,d_{\ell-1},d_{\ell}).$
\item[(ii)]
$(\A',m')$ is free with $\exp(\A',m')=(d_1,\ldots,d_{\ell-1},d_{\ell}-1).$
\item[(iii)]
$(\A'',m^{*})$ is free with $\exp(\A'',m^*)=(d_1,\ldots,d_{\ell-1}).$
\end{itemize}
In particular, 
if $(\A,m)$ and $(\A',m')$ are both free, then all the statements 
(i), (ii) and (iii) above hold.
\label{addition-deletion}
\end{theorem}

In general, the computation of Euler multiplicities $m^*$ is difficult without using 
a computer program. However, under some special condition, we can obtain $m^*$ in the 
following manner:

\begin{prop}[\cite{ATW2}, Proposition 4.1]\label{combemult}
Let $(\A,m)$ be a multiarrangement, 
$H_0 \in \A$ and $(\A'',m^*)$ the Euler restriction of 
$(\A,m)$ with respect to $H_0$. Let $X\in \mathcal{A}''$ and put $m_0=m(H_0)$. 
Suppose $k=|\mathcal{A}_X|$ and $m_1=\mathrm{max}\{ m(H)| H\in \mathcal{A}_X\backslash \{H_0\} \}$.

\begin{itemize}
\item[(1)] If $k=2$ then $m^*(X)=m_1$.

\item[(2)] If $2m_0\geq |m_X|$ then 
$m^*(X)=|m_X|-m_0$.

\item[(3)] If $2m_1\geq |m_X|-1$ then $m^*(X)=m_1$.

\item[(4)] If $|m_X|\leq 2k-1$ and $m_0>1$ then $m^*(X)=k-1.$

\item[(5)] If $|m_X| \leq 2k-2$ and $m_0=1$ then $m^*(X)=|m_X|-k+1.$

\item[(6)] If $m_X\equiv 2$ then $m^*(X)=k$.

\item[(7)] If $k=3$, $2m_0\leq |m_X|$, and $2m_1\leq |m_X|$ then $m^*(X) 
=\left\lfloor \frac{|m_X|}{2}\right\rfloor$.

\end{itemize}
\end{prop}

%For the proofs, details of these theorems, 
%especially the definitions of $(\A',m')$ and $(\A'',m^*)$, 
%see \cite{ATW2}. 
%Note that, for integers $d_1,\ldots,d_\ell$, $GMP_k(d_1,\ldots,d_\ell)$ 
%denotes the coefficient of $t^{\ell-k}$ of the polynomial 
%$\prod_{i=1}^\ell (t-d_i)$. 
Also,  
%in the application, we will 
to show the freeness of 
some deformations of the Coxeter arrangement, 
the following theorems by Ziegler in \cite{Z} and 
Yoshinaga in \cite{Y2} play central roles (see Section five). To 
introduce these results, let us review some definitions. 
Let $\A$ be a non-empty hyperplane arrangement and $H_0 \in \A$. 
The \textit{intersection lattice} $L(\A)$ of $\A$ is defined by 
$$
L(\A):=\{ \bigcap_{H \in \B} H|\B \subset \A\} 
$$
with the reverse inclusion as 
the partial ordering. 
For $X \in L(\A)$ the subarrangement $\A_X \subset \A$ is defined as the set 
$\{H \in \A|X \subset H\}$. $\A'$ is the \textit{deletion} of $\A$ 
with respect to $H_0$, defined by $\A':=\A \setminus \{H_0\}$. 
Also, $\A''$ is the \textit{restriction} of $\A$ with respect to $H_0$, 
defined by $\A'':=\{H' \cap H_0|H' \in \A'\}$. 
For each $X \in \A''$ we can associate the \textit{Ziegler multiplicity} $m_{H_0}$, 
defined in 
\cite{Z}, by $m_{H_0}(X):=|\{H' \in \A'|H' \cap H_0=X\}|$, and we call 
$(\A'',m_{H_0})$ the \textit{Ziegler restriction} with respect to $H_0$. 

\begin{theorem}[\cite{Z}]
In the above notation, if $\A$ is free with $\exp(\A)=(1,d_2,\ldots,d_\ell)$, then 
$(\A'',m_{H_0})$ is free with $\exp(\A'',m_{H_0})=(d_2,\ldots,d_\ell)$. 
\label{Ziegler}
\end{theorem}

\begin{theorem}[\cite{Y2}, Theorem 2.2]
In the above notation, assume that $\ell \ge 4$. Then 
$\A$ is free if and only if $(\A'',m_{H_0})$ is free and 
$\A_X$ is free for all $X \in L(\A'') \setminus \{\bigcap_{H \in \A} H\}$.
\label{Yoshinagamulti}
\end{theorem}

Next we introduce a criterion to check the non-freeness of 
multiarrangements, see \cite{ATW1} for the notation and details.

\begin{theorem}[\cite{ATW1}, Corollary 4.6]
If a multiarrangement $(\A,m)$ is free, then 
$GMP(k)=LMP(k)\ (1 \le k \le \ell)$, where $GMP(k)$ is the $k$-th global mixed product of 
$(\A,m)$ and 
$LMP(k)$ is the $k$-th local mixed product of $(\A,m)$.
\label{gmplmp}
\end{theorem}

The next proposition is useful to determine the non-freeness of multiarrangements, and 
the proof is the same as that for simple arrangements, 
see Theorem 4.37 in \cite{OT} for example.

\begin{prop}[\cite{A12}, Lemma 3.8] 
Let $(\A,m)$ be a multiarrangement and $X \in L(\A)$. 
%and 
%$m_X$ be a restriction of $m$ on $\A_X$. 
If $(\A,m)$ is free, then 
so is $(\A_X,m_X)$.
\label{localization}
\end{prop}

Next let us review the theory of a graphic arrangement and 
chordal graph by Stanley in \cite{S}. First, let us consider 
a subarrangement $\B$ of the Coxeter arrangement of type 
$A_\ell$. Then $\B$ can be uniquely characterized by using the graph $G$ consisting of 
the vertex set $V_G=\{1,2,\ldots,\ell+1\}$ and the 
set of non-directed edges $E_G$ in the following manner:

\begin{define}
For a graph $G$ as above, a \textit{graphic arrangement} $\A_G$ associated to the graph $G$ 
is defined by
$$
\A_G:=\{H_{ij}|\{i,j\} \in E_G\}.
$$
\end{define}

It is a natural problem to consider whether we can characterize the freeness of 
graphic arrangements in terms of the combinatorics of $G$. For that purpose, let us 
introduce the following graph.

\begin{define}
Let $G$ be a graph as above. 
A subgraph $C \subset G$ is a \textit{cycle} if $C$ consists of 
vertices $i_1,\ldots,i_s\ (s \ge 3)$ and $\{i_1,i_2\},\{i_2,i_3\},\ldots,
\{i_{s-1},i_s\},\{i_s,i_1\}$ are edges of $C$. A \textit{chord} of a cycle 
$C$ is an edge $\{i,j\}$ for non-consecutive vertices 
$i,j$ on the cycle $C$. 
% \in \{i_1,\ldots,i_s\}$. 
A graph $G$ is \textit{chordal} if every cycle $C \subset G$ with $|C| >3$ 
has a chord.
\label{chordalgraph}
\end{define} 

It is known that a graph is chordal if and only if 
its vertex set admits a vertex elimination order, see \cite{FG}. 
By using chordal graphs, Stanley gave a complete 
classification of free graphic arrangements as follows:

\begin{theorem}[\cite{S}]
A graphic arrangement $\A_G$ is free if and only if 
$G$ is chordal.
\label{Stanley}
\end{theorem}

For the rest of this article we give a 
generalization of Definition \ref{chordalgraph} and Theorem \ref{Stanley}.

\section{Signed-eliminable graphs}

In this section we introduce the theory of signed-eliminable graphs and give 
fundamental properties. 
This is a generalization of a chordal graph from the viewpoint of 
its vertex elimination ordering property. Recall 
the definition of a signed-eliminable graph in Definition \ref{bichordal} 
for the multi-braid arrangement. 
In the rest of this section, we introduce the theory of signed-eliminable 
graphs under the following setting. 

Let $G$ be a graph consisting of the vertex set $V_G$ with $|V_G|=\ell$ and 
the set of edges $E_G$ which has the decomposition $E_G=E_G^+ \cup E_G^-$ such that 
$E_G^+ \cap E_G^-= \emptyset$. 
For a subset $S \subset V_G,\ G|_S$ is the induced subgraph of 
$G$ with $V_{G|_S}=S$. 
%A map 
%$\nu:V_G \rightarrow \{1,2,\ldots,\ell\}$ is a 
%\textit{signed-elimination ordering} if 
%for $i,j,k \in V_G$ with $\nu(i), \nu(j)<\nu(k)$, 
%$G|_{\{i,j,k\}}$ is neither (1) nor (2) in Definition \ref{bichordal}.
We often consider that a sign $\sigma \in \{+,-\}$ is associated
to each edge in $E_G^{\sigma}$.
A signed graph $G$ is signed-eliminable if $V_G$ admits 
a signed-elimination ordering $\nu$. 

\begin{example}\label{ex:elimorder}
Let us classify all the signed-eliminable and non-signed-eliminable 
graphs with four vertices. Note that, by definition, 
the property that a graph is signed-eliminable is preserved even if we exchange the signs $+$ and $-$. Now 
the following graphs are signed-eliminable, where 
the numberings of vertices in the figure signify the corresponding signed-elimination 
ordering (we agree that 
an edge drawn in a single line belongs to $E_G^{\sigma}$ and that in a double line to 
$E_G^{-\sigma}\ (\sigma \in \{+,-\})$):
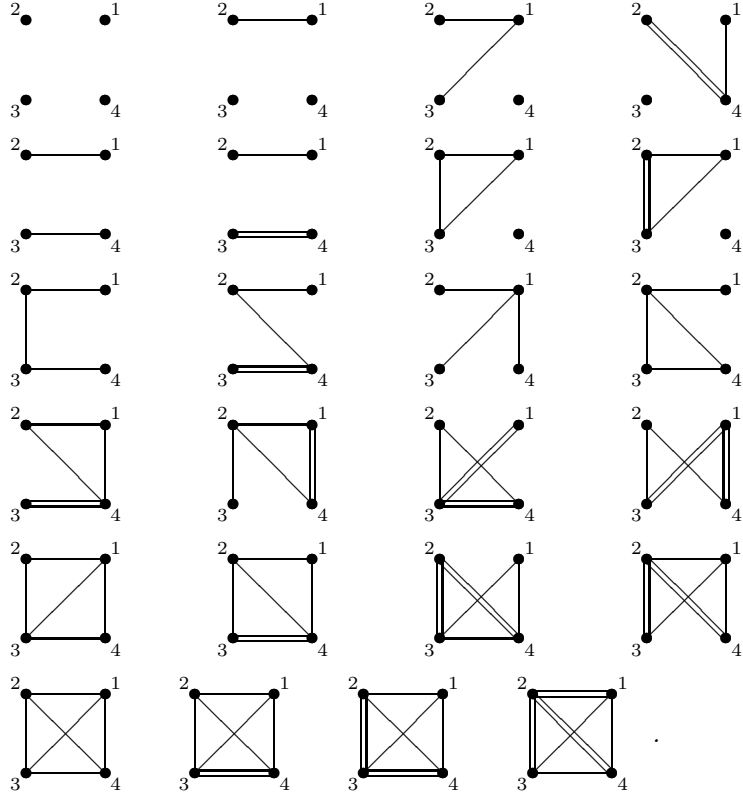
\medskip

\begin{center}
\parbox{0.8\textwidth}{
\GraphA{}{}{}{}{}{}
\GraphA{+}{}{}{}{}{}
\GraphA{+}{+}{}{}{}{}
\GraphA{}{}{+}{}{-}{}
\GraphA{+}{}{}{}{}{+}
\GraphA{+}{}{}{}{}{-}
\GraphA{+}{+}{}{+}{}{}
\GraphA{+}{+}{}{-}{}{}
\GraphA{+}{}{}{+}{}{+}
\GraphA{+}{}{}{}{+}{-}
\GraphA{+}{+}{+}{}{}{}
\GraphA{+}{}{}{+}{+}{+}
\GraphA{+}{}{+}{}{+}{-}
\GraphA{+}{}{-}{+}{+}{}
\GraphA{}{-}{}{+}{+}{-}
\GraphA{}{-}{-}{+}{+}{}
\GraphA{+}{+}{+}{+}{}{+}
\GraphA{+}{}{+}{+}{+}{-}
\GraphA{}{+}{+}{-}{-}{+}
\GraphA{+}{+}{+}{-}{-}{}
\GraphA{+}{+}{+}{+}{+}{+}
\GraphA{+}{+}{+}{+}{+}{-}
\GraphA{+}{+}{+}{-}{+}{-}
\GraphA{-}{+}{+}{-}{-}{+}.
}
\end{center}
\medskip

\break

The following graphs are not signed-eliminable:
\medskip

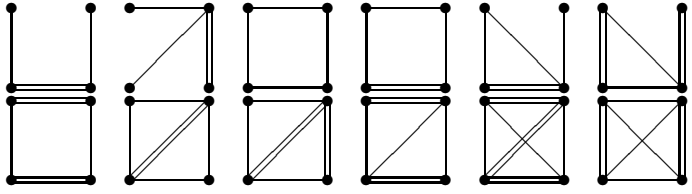
\begin{center}
\parbox{0.8\textwidth}{
\graphA{}{}{+}{+}{}{-}
\graphA{+}{+}{-}{}{}{}
\graphA{+}{}{+}{+}{}{+}
\graphA{+}{}{+}{+}{}{-}
\graphA{}{}{+}{+}{+}{-}
\graphA{}{}{-}{-}{+}{+}\\
\graphA{-}{}{+}{+}{}{-}
\graphA{+}{-}{+}{+}{}{+}
\graphA{+}{-}{-}{+}{}{+}
\graphA{-}{+}{+}{+}{}{-}
\graphA{-}{-}{+}{+}{+}{-}
\graphA{+}{+}{-}{-}{+}{+}.
}\end{center}
\end{example}

\begin{define}
Let $\nu$ be a signed-elimination ordering on $G$.
We define a \textit{$k$-th signed-eliminable filtration} of $G$ as 
a sequence of graphs 
$G_0,\ldots, G_m$
such that 

\begin{itemize}
\item $G_0=G|_{\{\nu^{-1}(1),\ldots,\nu^{-1}(k-1)\}}$,
\item $G_m=G|_{\{\nu^{-1}(1),\ldots,\nu^{-1}(k)\}}$,
\item $G_i$ is a subgraph of $G_{i+1}$ with decomposition of edge set induced by that of $G_{i+1}$,
\item $|E_{G_{i+1}} \setminus E_{G_{i}}|=1$, and 
\item $\nu|_{G_i}$ is a signed-elimination ordering on $G_i$ for each $i$.
\end{itemize}

For a signed-eliminable graph with $\ell$ vertices, 
we define a \textit{complete signed-eliminable filtration} of $G$ as 
a sequence of graphs 
$G_0,\ldots, G_m$
such that 
$G_{n_k},\ldots, G_{n_{k+1}}$
is  a $k$-th signed-eliminable filtration of $G$
for some $0=n_1 \leq n_2 \leq \cdots \leq n_{\ell+1}=m$.
\end{define}

\begin{rem}
The definition of a signed-eliminable graph with a signed-elimination ordering is 
just a generalization of the vertex elimination order  
on a non-signed graph.
%It is known that 
%a graph is chordal if and only if it has a vertex elimination order 
%(e.g., see \cite{FG}). 
Hence, 
from the viewpoint of Definition \ref{chordalgraph} and Theorem \ref{Stanley}, 
a signed-eliminable graph 
can be regarded as a generalization of a chordal graph. 
Theorem \ref{Nuida} in Section three also supports 
this generalization.
%that a signed-eliminable graph is a generalization of 
%a chordal graph from the viewpoint of the topology of graphs.
\label{remark}
\end{rem}

Let us investigate the properties of 
signed-eliminable graphs. The next proposition follows immediately by definition.

\begin{prop}
If some induced subgraph of $G$ is not signed-eliminable,
then $G$ is not signed-eliminable either.
\label{localelim}
\end{prop}

Now let us state the main theorem in this section, 
% to consider 
%the signed-eliminable graph. 
which will play the key role 
to characterize free multiplicities on the braid arrangement. 

\begin{theorem}
If $G$ is signed-eliminable, then $G$ always has a complete signed-eliminable filtration.
\label{key1}
\end{theorem}

%\begin{theorem}\label{lemma:welldefofinduction}
%Let $\nu$ be a signed-elimination order on $G$
%and $G_0,\ldots,G_i$ a signed-eliminable filtration on $G$.
%If $E_{G_{n+1}}=E_{G_{n}} \cup \{\{j_n,k_n\}\}$ and $\nu(j_n)<\nu(k_n)$,
%then for every $i \in V_G$ with $\nu(i)<\nu(k_n)$, 
%no induced subgraph $G_{n}|_{\{i,j_n,k_n\}}$ satisfies the following:
%\begin{itemize}
%\item $\{i,j_n\}$ and  $\{j_n,k_n\}$ are edges of the same color,
%and $\{i,k_n\}$ is an edge of another color.   
%\end{itemize}
%\end{theorem}

Roughly speaking, 
Theorem \ref{key1} ensures that we can 
always give an order on edges of 
a signed-eliminable graph 
which enables Addition-Deletion Theorem \ref{addition-deletion} 
work well. 
In the rest of this section we prove Theorem \ref{key1}.
% and 
%\ref{lemma:welldefofinduction}. 
For that purpose, we fix the following notation 
only in the rest of this section. 
Let $G$ be a signed-eliminable graph with $\ell$ vertices,
$\nu$ a signed-elimination ordering on $G$,
and $l \in V_G$ the vertex $\nu^{-1}(\ell)$.

\begin{lemma}\label{porder}
%Let $G$ be a bieliminable graph of order $\ell$,
%$\nu$ a bielimination order on $G$,
%and $l$ the vertex $\nu^{-1}(\ell)$. 
For $i,j \in V_G$, 
define the relation $i\prec j$ if
$\{i,j\}$ and $\{i,l\}$ are edges of the same sign 
and $\{j,l\}$ is an edge of the other sign.
Then the relation $\prec$ induces a partial order
on $\{i|\{i,l\}\in E_G\}$.
\end{lemma}

\noindent
\textbf{Proof}. 
First, let us show that $i_1 \prec i_2 \prec
i_3 \prec i_4$ implies $i_1 \prec i_4\ (i_a \in V_G)$. 
By symmetry, we may assume that 
$\{i_1,l\}, \{i_1,i_2\} \in E_G^+$ and 
$\{i_2,l\} \in E_G^-$. Then 
$\{i_2,i_3\} \in E_G^-,\ 
\{i_3,l\},\{i_3,i_4\} \in E_G^+$ 
and $\{i_4,l\} \in E_G^-$ by definition of 
$\prec$. Now if $\{i_1,i_4\} \not \in E_G^+$, then 
Example \ref{ex:elimorder} shows that 
$G|_{\{i_1,i_2,i_3,i_4\}}$ is not signed-eliminable, which contradicts 
Proposition \ref{localelim}. Hence $\{i_1,i_4\} \in E_G^+$, and $i_1 \prec i_4$. 

Now it suffices to show that there are no vertices 
$i_1,\ldots,i_n\ (n \ge 2)$ such that 
$i_1 \prec i_2 \prec \cdots \prec i_n \prec i_1$. 
If such vertices exist, then 
repeated use of the argument above implies that 
$i_1 \prec i_n \prec i_1$ (when $n$ is even) or 
$i_1 \prec i_2 \prec i_n \prec i_1$ (when $n$ is odd). 
However, this is impossible by definition of $\prec$. 
%A partially ordering is well-defined
%if there are no vertices $i_1,\ldots,i_n $ such that
%$i_1\prec i_2\prec \cdots\prec i_n\prec i_1\ (n \ge 2)$. 
%The definition of $\prec$ implies that $n$ is even. If $n=2$ then 
%this contradicts to the definition of a signed-eliminable graph. Hence 
%assume that $n \ge 4$. Then 
%this order 
%of vertices can occur only when they satisfy the conditions 
%in Lemma \ref{lemma:nocycle}, which completes the proof.
\owari

\begin{lemma}\label{key-1}
%Let $G$ be a bieliminable graph of order $\ell$,
%$\nu$ a bielimination order on $G$,
%$l$ the vertex $\nu^{-1}(\ell)$, 
Let $j$ be a maximal vertex of the poset $\{i|\{i,l\}\in E_G\}$
defined by $\prec$ in Lemma \ref{porder} and 
$G'$ the graph obtained from $G$
by deleting the edge $\{j,l\}$.
Then $G'$ is also signed-eliminable with the same signed-elimination ordering 
$\nu$.
\end{lemma}

\noindent
\textbf{Proof}. 
By the definition of the signed-eliminable graph, 
it is sufficient to consider 
the induced subgraph $G'|_{\{i,j,l\}}$ for any $i$ with 
$\nu(i) <\nu(l)$. The  
classification of every possible case for $G|_{\{i,j,l\}}$ 
shows that the induced subgraph $G'|_{\{i,j,l\}}$ does not satisfy 
the conditions of 
Definition \ref{bichordal}
%not signed-eliminable
only if 
$\{i,j\}$ and $\{j,l\}$ are edges of the same sign
and $\{i,l\}$ is an edge of the other sign in $G|_{\{i,j,l\}}$.
%\footnote{
%�����A�O�̓��ӂ�$(i,j)(j,l)$�Ō��낪
%$(i,l)$�ł͂Ȃ��ł��傤���H�@�m�F���肢���܂�}.
However, we 
have assumed that $j$ is a maximal vertex of the poset $\{i|\{i,l\}\in E_G\}$
defined by $\prec$, which completes the proof.
\owari
\medskip

\noindent
\textbf{Proof of Theorem \ref{key1}}. 
Apply Lemma \ref{key-1} repeatedly 
to edges $\{\{i,l\} \in E_G|\nu(i)<\nu(l)\}$. \owari
\medskip

%\noindent
%\textbf{Proof of Theorem \ref{lemma:welldefofinduction}}. 
%The following lemma completes the proof of Theorem \ref{lemma:welldefofinduction}. \owari

%\begin{lemma}
%Let $G$ and $G'$ be signed-eliminable graphs of order $\ell$,
%$\nu$ a signed-elimination order on $G$ and $G'$, 
%and $l$ the vertex $\nu^{-1}(\ell)$. 
%If $E_{G'}=\{\{j,l\}\} \cup E_{G}$,
%then for $i \in V_G$ with $\nu(i) < \nu(l)$, 
%no induced subgraph $G_{n}|_{\{i,j,l\}}$ is

%\begin{itemize}
%\item $\{i,j\}$ and  $\{j,l\}$ are edges of the same color,
%and $\{i,l\}$ is an edge of another color.   
%\end{itemize}
%\end{lemma}

%\noindent
%\textbf{Proof}. 
%Let $\{i,j\}$ and  $\{j,l\}$ 
%be edges of the same color in $G$,
%and $\{i,l\}$  an edge of another color  in $G$.  
%Then 
%the induced subgraph $G'|_{\{i,j,l\}}$ is

%\begin{itemize}
%\item $\{i,j\}$ and $\{i,l\}$ are edges of different
%      colors,
% and $\{j,l\}$ is disjoint.
%\end{itemize}

%Recalling that  $G'$ is signed-eliminable with the same ordering $\nu$, this is contradiction.
%\owari

\section{Characterization of signed-eliminable graphs}
In this section we quote a characterization of signed-eliminable graphs from \cite{N}. To state it, 
let us introduce the following two definitions.

\begin{define}[\cite{N}, Definition 4.4]
Let $G$ be a graph with the set of vertex $V_G$ and two sets of edges 
$E_G^+$ and $E_G^-$ as in the previous section, and $\sigma \in \{+,-\}$.
\begin{itemize}
\item[(1)] 
A sequence $(v_1,v_2,\ldots,v_n;\omega)\ (n \ge 3)$ of vertices in $G$ 
is a $(\sigma$-)\textit{mountain} if 
$\{v_i,v_{i+1}\} \in E_G^{-\sigma}$ for $1 \le i \le n-1$, 
$\{\omega,v_i\} \in E_G^\sigma$ for $2 \le i \le n-1$ and any other 
pair of vertices is not joined by an edge. 

\item[(2)] 
A sequence $(v_1,v_2,\ldots,v_n;\omega_1,\omega_2)\ (n \ge 2)$ of vertices in $G$ 
is a $(\sigma$-)\textit{hill} if 
$\{v_i,v_{i+1}\} \in E_G^{-\sigma}$ for $1 \le i \le n-1$, 
$\{\omega_1,\omega_2\} \in E_G^\sigma$, 
$\{\omega_1,v_i\} \in E_G^\sigma$ for $1 \le i \le n-1$, 
$\{\omega_2,v_i\} \in E_G^\sigma$ for $2 \le i \le n$ and any other 
pair of vertices is not joined by an edge. 
\end{itemize}
\label{mountainhill}
\end{define} 

By using chordality, mountains, hills, and Example \ref{ex:elimorder}, 
a characterization of signed-eliminable graphs is given as follows.

\begin{theorem}[\cite{N}, Theorem 5.1]
Let $G$ be a signed graph. 
Then $G$ is signed-eliminable if and only if the following three conditions are 
satisfied:
\begin{itemize}
\item[(C1)] Both graphs $(V_G,E_G^+)$ and $(V_G,E_G^-)$ are chordal.
\item[(C2)] Any induced subgraph of $G$ with four vertices is 
signed-eliminable.
\item[(C3)] $G$ contains no mountains nor hills.
\end{itemize}
\label{Nuida}
\end{theorem}

For details of Theorem \ref{Nuida}, see \cite{N}. Theorem \ref{Nuida} 
plays the key role for the proof of the ``only if'' part of Theorem \ref{freemulti}. 
Note that, if $E_G^-=\emptyset$, then Theorem \ref{Nuida} asserts the 
well-known equivalence between a chordal graph and a graph with a vertex elimination ordering.

\section{Proof of Theorem \ref{freemulti}}
In this section we apply the theory of signed-eliminable graphs 
to prove Theorem \ref{freemulti}. Since the proof is the same, 
we only prove the case when the condition (a) in Theorem \ref{freemulti} is 
satisfied. 

First, let us prove the ``if'' part. 
Let $G$ be a signed-eliminable graph with a signed-elimination 
ordering $\nu:V_G \rightarrow \{1,2,\ldots,\ell+1\}$. 
By an appropriate change of coordinates, we may assume that 
$\nu(v_i)=i$ for all $i$. Then let us identify $v_i$ with $i$ for all $i$ in this 
proof. 
Hence the order of vertices $V_G=\{1,2,\ldots,\ell+1\}$ is already a 
signed-elimination ordering. 
When $E_G=\emptyset$, the theorem %is proved in 
%\cite{AN2}, or 
can be proved 
by using the argument below with the signed-eliminable graph $G$ consisting of 
$V_G=\{1,2,\ldots,\ell+1\}$ and $E_G=E_G^+=\{\{i,j\}|j=1,\ldots,\ell+1,j \neq i\}$ for 
a fixed $i$. We prove the statement 
by induction on $\ell$. When $\ell=1$ there is nothing to prove. 
If $\ell=2$ then the result in \cite{Waka} completes 
the proof. Assume that $\ell >2$. Also, assume that 
$\A_\ell(n_1,\ldots,n_{\ell+1})[G|_{\{1,2,\ldots,s-1\}}]$ is free with 
exponents $(0,N+\widetilde{\deg}_2,\ldots,N+\widetilde{\deg}_{s-1},N,\ldots,N)$ for 
some $s,\ 2 \le s \le \ell+1$. 
By Theorem \ref{key1}, there exists an $s$-th filtration 
$G_0^s,\ldots,G_{f(s)}^s$ of 
$G$ with $E_{G_{i+1}^s} \setminus E_{G_i^s}=\{\{s, j_{i}\}\}\ (j_i<s)$. 
Consider the Euler restriction $(\A'',m^*)$ of the multiarrangement 
$\A_\ell(n_1,\ldots,n_{\ell+1})[G_{i+1}^s]$ onto the hyperplane 
$H_{s j_i}\ (i=0,1,\ldots,f(s)-1)$. 
Combining Theorem \ref{addition-deletion} and Proposition \ref{combemult} with 
Definition \ref{bichordal} and Theorem \ref{key1}, 
% and \ref{lemma:welldefofinduction}, 
the lemma below follows 
immediately.

\begin{lemma}\label{preserve}
In the notation above, let $t\in V_G$ with $t<s$. 
If $(\A'',m^*)$ is the Euler restriction with respect to 
$H_{sj_i}$, then 
$$
m^*(H_{tj_i})=m^*(H_{ts})=
3k+n_{j_i}+n_s+n_t+m_G(H_{tj_i}).
$$
\end{lemma}

Then Lemma \ref{preserve} implies that the Euler restriction 
$(\A'',m^*)$ is equal to the following multiarrangement:
$$
\A_{\ell-1}(n_1,\ldots,n_{j_i-1}, n_{j_i}+n_s+k,n_{j_i+1},\ldots,
n_{s-1},n_{s+1},\ldots,n_{\ell+1})[G|_{\{1,2,\ldots,s-1\}}]. 
$$
Proposition \ref{localelim} and Theorem \ref{key1} imply 
that $G|_{\{1,2,\ldots,s-1\}}$ is also signed-eliminable with 
a signed-elimination ordering $\{1,2,\ldots,s-1\}$. Hence 
the induction hypothesis shows that $(\A'',m^*)$ is free with exponents 
$(0,N+\widetilde{\deg}_2,\ldots,N+\widetilde{\deg}_{s-1},N,\ldots,N)$. Then 
Addition-Deletion Theorem \ref{addition-deletion} completes the proof of the 
``if'' part. 

Next we prove the ``only if'' part. Assume that 
$G$ is not signed-eliminable. Then   
Theorem \ref{Nuida} implies that  
$G$ does not satisfy the conditions (C1), (C2) or (C3). Also 
identify $v_i$ with $i$ for all $i$ in this proof. 
We will prove that $\A_\ell(n_1,\ldots,n_{\ell+1})[G]$ is not free in each of these three cases. 
To prove it, let us introduce a definition used only in this proof. 
A signed graph $G$ is \textit{free} if the associated multi-braid 
arrangement $\A_\ell(n_1,\ldots,n_{\ell+1})[G]$ is free. 
First, assume that 
$G$ does not satisfy the condition (C2). Then  
$G$ contains some non-signed-eliminable subgraph with 
four vertices. 
By Example \ref{ex:elimorder}, such a graph is one of the following:
% or 
%thier dual (i.e., the graph whose color of every edge is exchanged from $+$ to $-$ and 
%vice versa):
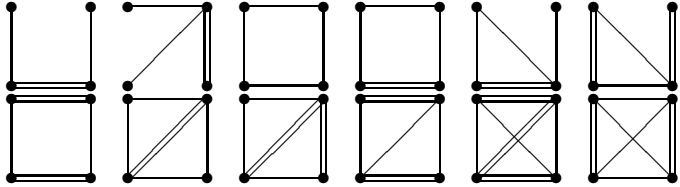
\begin{center}
\parbox{0.8\textwidth}{
\graphA{}{}{+}{+}{}{-}
\graphA{+}{+}{-}{}{}{}
\graphA{+}{}{+}{+}{}{+}
\graphA{+}{}{+}{+}{}{-}
\graphA{}{}{+}{+}{+}{-}
\graphA{}{}{-}{-}{+}{+}\\
\graphA{-}{}{+}{+}{}{-}
\graphA{+}{-}{+}{+}{}{+}
\graphA{+}{-}{-}{+}{}{+}
\graphA{-}{+}{+}{+}{}{-}
\graphA{-}{-}{+}{+}{+}{-}
\graphA{+}{+}{-}{-}{+}{+}.
}\end{center}
By Proposition 
\ref{localization} it suffices to show that 
these graphs are not free. For that purpose, we use two theorems, i.e., 
Theorems \ref{addition-deletion} and \ref{gmplmp}. First, prove the non-freeness of the graphs 
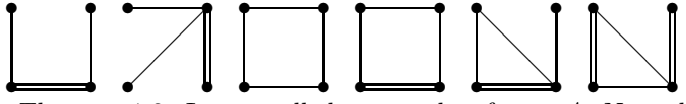
\begin{center}
\parbox{0.8\textwidth}{
\graphA{}{}{+}{+}{}{-}
\graphA{+}{+}{-}{}{}{}
\graphA{+}{}{+}{+}{}{+}
\graphA{+}{}{+}{+}{}{-}
%\graphA{-}{}{+}{+}{}{-}
\graphA{}{}{+}{+}{+}{-}
\graphA{}{}{-}{-}{+}{+}
%\graphA{+}{-}{+}{+}{}{+}
%\graphA{+}{-}{-}{+}{}{+}
%\graphA{-}{+}{+}{+}{}{-}
%\graphA{-}{-}{+}{+}{+}{-}
%\graphA{+}{+}{-}{-}{+}{+}.
}\end{center}
by using Theorem \ref{addition-deletion}. Let us call these graphs of type $A$. 
Note that, by deleting 
an appropriate edge from graphs of type $A$, we can obtain 
signed-eliminable graphs as follows:
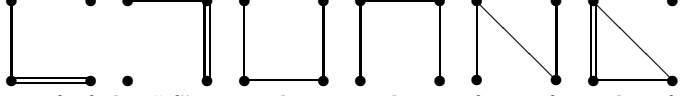
\begin{center}
\parbox{0.8\textwidth}{
\graphA{}{}{}{+}{}{-}
\graphA{+}{}{-}{}{}{}
\graphA{}{}{+}{+}{}{+}
\graphA{+}{}{+}{+}{}{}
%\graphA{-}{}{+}{+}{}{-}
\graphA{}{}{+}{+}{+}{}
\graphA{}{}{}{-}{+}{+}
%\graphA{+}{-}{+}{+}{}{+}
%\graphA{+}{-}{-}{+}{}{+}
%\graphA{-}{+}{+}{+}{}{-}
%\graphA{-}{-}{+}{+}{+}{-}
%\graphA{+}{+}{-}{-}{+}{+}.
}\end{center}
By the proof of the ``if'' part, these graphs are free. 
If graphs of type A are also free, then Theorem \ref{addition-deletion} 
implies that $\exp(\A'',m^*) \subset \exp(\A',m')$ as multisets, which contradicts  
the results in \cite{Waka}, Proposition \ref{combemult} and what is proved in the ``if'' part. Hence graphs of 
type $A$ are not free. 
%The same proof is valid for dual graphs of type $A$. 
Next let us prove the non-freeness of the remaining graphs 
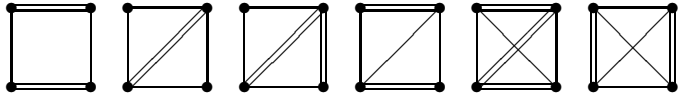
\begin{center}
\parbox{0.8\textwidth}{
%\graphA{}{}{}{+}{}{-}
%\graphA{+}{}{-}{}{}{}
%\graphA{}{}{+}{+}{}{+}
%\graphA{+}{}{+}{+}{}{}
\graphA{-}{}{+}{+}{}{-}
%\graphA{}{}{+}{+}{+}{}
%\graphA{}{}{}{-}{+}{+}
\graphA{+}{-}{+}{+}{}{+}
\graphA{+}{-}{-}{+}{}{+}
\graphA{-}{+}{+}{+}{}{-}
\graphA{-}{-}{+}{+}{+}{-}
\graphA{+}{+}{-}{-}{+}{+}
}\end{center}
by using Theorem \ref{gmplmp}. 
Let us call these graphs of type $B$ and give a name $B_1,B_2,\ldots,B_6$ to each 
of these graphs from the left. Assume that graphs of type $B$ are 
free. 
Also, assume that a single line edge corresponds to an edge in $E_G^+$ and 
a double line edge to that in $E_G^-$.  
Let $G_i$ (resp.$\ L_i)$ denote the 2nd global (resp.$\ $local) mixed product 
of $\A_3(n_1,n_2,n_3,n_4)[B_i]$. Then we can compute these values according to 
\cite{ATW1} as follows (where $N=\sum_{i=1}^4 n_i)$:
\begin{itemize}
\item[$B_1:$] $G_1 \le 48k^2 +24kN +3N^2<L_1=48k^2 +24kN +3N^2+2$.
\item[$B_2:$] $G_2 \le 48k^2 +24kN +3N^2+6N+24k+3<L_2=48k^2 +24kN +3N^2+6N+24k+4$.
\item[$B_3:$] $G_3 \le 48k^2 +24kN +3N^2+8k+2N<L_3=48k^2 +24kN +3N^2+8k+2N+1$.
\item[$B_4:$] $G_4 \le 48k^2 +24kN +3N^2+8k+2N<L_4=48k^2 +24kN +3N^2+8k+2N+2$.
\item[$B_5:$] $G_5 \le 48k^2 +24kN +3N^2<L_5=48k^2 +24kN +3N^2+1$.
\item[$B_6:$] $G_6 \le 48k^2 +24kN +3N^2+4N+16k+1<L_6=48k^2 +24kN +3N^2++4N+16k+3$.
\end{itemize}
Hence Theorem \ref{gmplmp} implies contradictions, which show that these 
graphs are not free. Since the same proof as the above is valid when 
the signs of single and double lines are exchanged, graphs of type $B$ are not free, 
which shows that 
every non-signed-eliminable graph with four 
vertices is not free. 
%Then the non-freeness of these submultiarrangements follows immediately by 
%Proposition \ref{localization} and case-by-case calculation through 
%the Local-Global formula Theorem \ref{gmplmp} for $2$nd global and local 
%mixed products (e.g., see \cite{ATW1}), and addition-deletion Theorems 
%\ref{addition-deletion}. \footnote{LMPGMP�������Ă��ǂ�.}

Next assume that the condition (C1) is not satisfied. Then 
there exists a subgraph $C \subset G$ such that 
$|C| \ge 4$ and $(V_C, E_G^{\sigma} \cap E_C)$ is a cycle without 
chords of the sign $\sigma \in \{+,-\}$. 
Because of the symmetry we may 
assume that $\sigma=+$. 
Moreover, Proposition \ref{localization} implies that 
it is sufficient to show that $C$ or its subgraph 
is not free. We prove the non-freeness 
by induction on $\ell \ge 2$. If $\ell=2$ then there is nothing 
to prove, so assume that $\ell >2$. 
If $|C|=4$ then Example \ref{ex:elimorder}
implies that $C$ is not signed-eliminable, hence 
the above arguments imply the non-freeness. Assume that 
$|C| > 4$. First, assume that  there are no chords in $E_C^+ \cup E_C^-$. 
When $|C| < \ell+1$, the induction hypothesis 
completes the proof. So we may assume that $|C|=\ell+1$. 
We may also assume that $\{\{1,2\},\{2,3\},\ldots, \{\ell,\ell+1\},\{\ell+1,1\}\} = E_C^+=E_C$. 
Define a subgraph $C' \subset C$ which is obtained from $C$ by deleting the edge 
$\{\ell+1,1\}$. Note that $C'$ is signed-eliminable. 
Then the ``if'' part of Theorem \ref{freemulti} implies that 
$\A_\ell(n_1,\ldots,n_{\ell+1})[C']$ is free with exponents 
$(0,N+1,\ldots,N+1)$. If $\A_\ell(n_1,\ldots,n_{\ell+1})[C]$ is free, then 
every statement in Theorem \ref{addition-deletion} holds. Let us consider the 
Euler restriction of $\A_\ell(n_1,\ldots,n_{\ell+1})[C]$ onto $x_{\ell+1}-x_1=0$. 
Then the Euler restriction is equivalent to 
$
\A_{\ell-1}(n_1+n_{\ell+1}+k,n_2,\ldots,n_{\ell})[C''],
$
where $C''$ is a cycle with $V_{C''}=\{1,2,\ldots,\ell\}$ and 
$E_{C''}=E_{C''}^+=\{\{1,2\},\{2,3\},\ldots,\{\ell-1,\ell\},\{\ell,1\}\}$. 
%all of whose edges are colored 
%by $+$. 
If $\ell=3$, then \cite{Waka} implies the contradiction on the exponents.  
If $\ell >3$ then the induction hypothesis shows that the 
Euler restriction is not free, which is also a contradiction. 

So we may assume that the cycle $C$ contains a chord whose sign is $-$. 
Use the same notation in the above paragraph and assume that the chord is 
$\{i,j\}$, where $i$ and $j$ are non-consecutive vertices in $V_C$ 
with $i < j$. Also we may assume that 
$i \neq 1$ and $j \neq \ell+1$. 
Then we obtain two new graphs $C_1$ and $C_2$ as induced subgraphs of 
$C$ with $V_{C_1}=\{1,2,\ldots,i,j,j+1,\ldots,\ell+1\}$ and  
$V_{C_2}=\{i,i+1,\ldots,j\}$ 
respectively. 
If $|C_1|=4$ or $|C_2|=4$, then the previous argument for the non-freeness 
of non-signed-eliminable graphs with four vertices and Example \ref{ex:elimorder} complete the proof. 
If, for example, $|C_1| > 4$, then we may take a subgraph $C_1' \subset C_1$ whose vertices consist of 
$\{i-1,i,j,j+1\}$. If $E_{C_1'}=\{\{i-1,i\},\{i,j\},\{j,j+1\}\}$, 
then $C_1'$ is not signed-eliminable with 
four vertices, hence not free as we have already proved in the above. 
If there is some other edge in $C_1'$, then the assumption 
implies that edge has to be signed by $-$. If that edge is  
$\{i-1,j\}$, then consider the induced subgraph $C_{11} \subset C_1$ whose vertices consist of 
$\{1,2,\ldots,i-1,j,j+1,\ldots,\ell+1\}$ and apply the same arguments above. Then finally, we 
obtain a non-signed-eliminable, hence non-free subgraph with four vertices, which completes the proof. 

Finally, assume that the condition (C3) is not satisfied. 
Because the proof is the same, let us assume that $G$ contains 
a $(+)$-mountain $C=(v_1,v_2,\ldots,v_s;\omega) \subset G\ (s \ge 3)$. 
By Proposition \ref{localization} it suffices to show that 
$C$ is not free.
% which consists of a ``top'' vertex 
%$\omega$ and ``flat'' vertices $v_1,v_2,\ldots,v_s\ $ (see 
%\cite{N} for details). 
If $s=3$, then Example \ref{ex:elimorder} implies that 
$C$ is not signed-eliminable. Hence the first argument of the ``only if'' part 
of Theorem \ref{freemulti} shows the non-freeness. 
Assume that $s > 3$. Consider the subgraph $C' \subset C$ which is obtained from $C$ by 
deleting the vertex $v_s$ and the edge $\{v_{s-1},v_s\} \in E_G^-$. Then $C'$ 
has a signed-elimination ordering whose  
$k$-th filtration is given by first adding $\{w,v_{k-1}\}$ and second adding 
$\{v_{k-2}, v_{k-1}\}$, hence $C'$ is free by the ``if'' part of Theorem \ref{freemulti}. If 
$\A_\ell(n_{1},\ldots,n_{v_{\ell+1}})[C]$ is free, then  
Theorem \ref{addition-deletion} implies that the Euler restriction $(\A'',m^*)$
of $\A_\ell(n_1,\ldots,n_{\ell+1})[C]$ onto $H_{v_{s-1} v_s}$ is also free. 
However, Proposition \ref{combemult} implies that the Euler restriction 
$(\A'',m^*)$ 
corresponds to the graph of the mountain 
$(v_1,v_2,\ldots,v_{s-1};\omega)$, 
%whose top vertex is $w$ and flat vertices are $\{v_1,\ldots,v_{n-1}\}$, 
hence not free by the induction hypothesis.

When $G$ contains a hill, the same proof as the 
above can be applied, which completes the proof of 
Theorem \ref{freemulti}.
\owari
\medskip

Since exponents do not depend on a choice of a basis as the multiset, the next corollary 
follows immediately from Theorem \ref{freemulti}.

\begin{cor}
If
$G$ is signed-eliminable, then 
$\widetilde{\deg}_1=0$ and $(\widetilde{\deg}_1,\widetilde{\deg}_2,\ldots,\widetilde{\deg}_{\ell+1})$ does not depend on 
a choice of a signed-elimination ordering as the multiset. 
\label{indepofdeg}
\end{cor}

In \cite{ATW1}, a characteristic polynomial $\chi(\A,m,t)$ of multiarrangements is defined and 
the factorization theorem is proved. In general, the computation of $\chi(\A,m,t)$ is difficult, but 
if $(\A,m)$ is free, then we can easily compute it by the factorization. So when $G$ is signed-eliminable, 
we can calculate its characteristic polynomial as follows:

\begin{cor}
Let $(\A,m)=\A_\ell(n_1,\ldots,n_{\ell+1})[G]$ be the same as in 
Theorem \ref{freemulti}. 
Define $(\A,\tilde{m})$ by $\tilde{m}(H_{ij}):=2k+n_i+n_j-m_G(H_{ij})$. 
\begin{itemize}
\item[(1)] 
Let $k>0$. Then $(\A,m)$ is free if and only if $(\A,\tilde{m})$ is free.
\item[(2)] 
If $G$ is signed-eliminable, then 
$$
\chi(\A,m)=t\prod_{i=2}^{\ell+1} (t-N-\widetilde{\deg}_{i})
$$
and
$$
\chi(\A,\tilde{m})=t\prod_{i=2}^{\ell+1} (t-N+\widetilde{\deg}_{i}).
$$
\end{itemize}
\label{cha}
\end{cor}

Corollary \ref{cha} shows that there exists a duality of exponents 
of free multi-braid arrangements as mentioned in \cite{AY2}.

\section{Conjecture of Athanasiadis}
In this section we apply the results in previous sections to a conjecture 
of Athanasiadis in \cite{Ath2}. To state it, let us introduce 
some notation. 

Let us consider an affine arrangement in $V^{\ell+1}$ defined by 
\begin{eqnarray}\label{defath}
x_i-x_j=-k-\epsilon(i,j),-k,-(k-1),\ldots,k,k+\epsilon(j,i)  \\
 (1 \le i<j \le \ell+1), \nonumber
\end{eqnarray}
where $k \in \Z_{\ge 0}$ and $\epsilon(i,j)=0$ or $1$. Note that 
in this section, we distinguish $(i,j)$ and $(j,i)$ as explained later. 
Such arrangements are examples of 
\textit{deformations of the braid arrangement}, a class of arrangements 
first investigated systematically by Stanley in \cite{S2}. 
From the viewpoint of the combinatorics and freeness, these 
arrangements 
have been extensively studied by Athanasiadis \cite{Ath0}, 
\cite{Ath1}, \cite{Ath2}, Edelman and Reiner \cite{ER2}, 
Postnikov and Stanley \cite{PS}, Yoshinaga \cite{Y2} and many other authors. 
The main focus of these authors is on the characteristic polynomial of these arrangements. 
%One of the main researches  
%for them is to calculate their characteristic polynomial, and 
Because of Terao's factorization theorem, it is important 
to consider the freeness of these arrangements. 

Now let us go back to the deformation (\ref{defath}).
A useful way to consider this arrangement is introduced by Athanasiadis in 
\cite{Ath0}. Consider the directed graph $G$ consisting of the vertex set 
$V_G=\{1,2,\ldots,\ell+1\}$ and the set of directed edges 
$E_G \subset \{(i,j)|1 \le i,\ j \le \ell+1\}$. Here the edge 
$(i,j)$ is the arrow from $i$ to $j$. If we define 
\[
\epsilon(i,j):=
\left\{
\begin{array}{rl}
1 & \mbox{if}\ (i,j) \in E_G ,\\
0 & \mbox{if}\ (i,j) \not \in E_G, 
\end{array}
\right.
\]
%$\epsilon(i,j)$ to be zero if there are no edges between $i$ and $j$, and to be 
%one if there is an arrow from $i$ to $j$, 
then every affine arrangement 
above can be expressed by these directed graphs. For such a graph $G$ let $\A_G$ denote 
the corresponding arrangement of the form (\ref{defath}). In \cite{Ath0}, Athanasiadis gave 
a splitting formula of the characteristic polynomial of $\A_G$ when $G$ satisfies 
the following two conditions: 

\begin{itemize}
\item[(A1)] 
For every triple $i,j,h$ with $i,j<h$, it holds that, if $(i,j) \in E_G$, then 
$(i,h) \in E_G$ or $(h,j) \in E_G$.
\item[(A2)] For every triple $i,j,h$ with $i,j<h$, it holds that, 
if $(i,h) \in E_G$ and $(h,j) \in E_G$ then $(i,j) \in E_G$. 
\end{itemize} 

Athanasiadis also gave the following conjecture. 

\begin{conjecture}[\cite{Ath2}, Conjecture 6.6]\label{conjath}
Let $k=0$ in the deformation 
(\ref{defath}). Then the coning 
$c\A_G$ of $\A_G$ is free if and only if $G$ satisfies conditions (A1) and (A2).
\end{conjecture}

%In fact, Conjecture \ref{conjath} is a slight generalized version of 
%that in \cite{Ath2}, in which Athanasiadis gave the same statement when 
%$k=0$. However, as seen in Theorem \ref{freemulti} or \cite{AY2}, 
%the shifting by $2k$ seems to give no changes on the freeness of 
%these deformations. So we 
%state a conjecture of Athanasiadis in the form above. 
In the rest of this section let us prove that 
(A1) and (A2) are sufficient conditions in Conjecture \ref{conjath} 
in more general setting. First, let us prove the following.

\begin{prop}
Let $H_{\infty} \in c \A_G$ be the infinity hyperplane 
of the coning $c\A_G$ of $\A_G$ in (\ref{defath}). If $G$ satisfies (A1) and (A2), then 
the Ziegler restriction $(\A'', m_{H_{\infty}})$ with respect to 
$H_{\infty}$ 
is of the form $\A_{\ell}(n_1,\ldots,n_{\ell+1})[G']$ 
for some $n_1,\ldots,n_{\ell+1}$ and 
signed-eliminable graph $G'$. In particular, it is free.
\label{key2}
\end{prop}

\noindent
\textbf{Proof.} 
Note that the 
signed-eliminability is a local condition. In other words,  
that can be determined by checking the behavior of edges between every ordered 
triple of vertices $i,j<h$. 
Hence the proposition follows immediately by conditions (A1), (A2), the 
definition of a signed-eliminable graph and Theorem \ref{freemulti}. \owari
\medskip

\begin{theorem}
%The ``if'' part of Conjecture \ref{conjath} is true, i.e., 
In the deformation (\ref{defath}), 
$c\A_G$ is free if $G$ satisfies (A1) and (A2). 
In particular, the ``if'' part of Conjecture \ref{conjath} is true. 
\end{theorem}

\noindent
\textbf{Proof}.  
Induction on $\ell \ge 1$. When $\ell=1$, there is nothing to prove. 
If $\ell=2$ then 
the classification in \cite{A1} completes the proof. 
Assume that $\ell \ge 3$. 
By Theorem \ref{Yoshinagamulti} 
%Theorem \ref{key1} 
and Proposition \ref{key2}, 
it suffices to show that $(c \A_G)_X$ is free for any 
$X \in L(c\A_G)$ with $\bigcap_{H \in c \A_G} H \subsetneq X \subset H_{\infty}$. 
Again, recall that 
conditions (A1) and (A2) are local and note that 
$(c \A_G)_X$ decomposes into the direct product of 
the empty arrangement and the arrangement $c \A_{G'}$, where 
$G'$ is some directed graph. 
In fact, if $X=\{x_{i_1}=x_{i_2}=\cdots=x_{i_s}\} \cap H_{\infty}$, then 
$G'$ is the induced subgraph of $G$ with $V_{G'}=\{i_1,\ldots,i_s\}$. 
Then again the locality of (A1) and (A2) implies that  
$G'$ also satisfies conditions (A1) and (A2). Since 
$\rank (c \A_{G'}) < \rank (c \A_{G})$, the 
induction hypothesis implies that $c \A_{G'}$ is free. Hence  
$(c \A_G)_X$ is also free, which completes the proof. \owari

% \vspace{5mm}

%\noindent
%Takuro Abe\\
%Department of Mathematics\\
%Hokkaido University\\
%Sapporo 060--0810, Japan\\
%abetaku@math.sci.hokudai.ac.jp
%
%\bigskip
%
%\noindent
%Koji Nuida\\
%Research Center for Information Security (RCIS)\\
%National Institute of Advanced Industrial Science and Technology (AIST)\\
%Tokyo 101--0021, Japan\\
%k.nuida@aist.go.jp
%
%%\noindent%
%
%%{\em e-mail address}\ : \ 
%
%\bigskip
%
%\noindent
%Yasuhide Numata\\
%Department of Mathematics\\
%Hokkaido University\\
%Sapporo 060--0810, Japan\\
%nu@math.sci.hokudai.ac.jp

\appendix
\section{Erratum of Theorem \ref{freemulti}	}

Theorem 0.3 in this paper is not correct as it was stated. Explicitly, the statement does not hold when 
the graph $G$ satisfies the condition (iii) in the paper. See the following example by Michael Dipasquale.
\medskip

\noindent
\textbf{Example}.\\ 
Consider the multiarrangement $(x_1-x_2)(x_1-x_3)(x_1-x_4)(x_2-x_3)^4(x_2-x_4)^2(x_3-x_4)^2$ on the $A_3$-type. 
This multiplicity corresponds to the case $k=0,\ n_1=0,\ n_2=n_3=2,\ n_4=1$ with the 
graph $G$ such that $E_G^+=\emptyset,\ E_G^{-}=\{\{1,2\},\ 
\{2,4\},\ \{3,4\},\ \{1,3\}\}$. Since $E_{G}^-$ is not signed-eliminable, this has to be not free. However, 
by Theorem 5.10 in \cite{ATWA}, this is a supersolvable multiarrangement, hence free with exponents $(3,4,4)$, which 
contradicts % Hence t
Theorem 0.3 (iii)  in the paper. 
%is not correct. 
\medskip

Under the terminology of the paper, the correct statement of Theorem 0.3 should be as follows:
\medskip

\noindent
\textbf{Theorem 0.3}\\
Let $\A$ be the braid arrangement in $V^{\ell+1}$, $G$ a signed graph and let 
$m_G$ be the map in Definition 0.1. Let $k,n_1,\ldots,n_{\ell+1}$ be non-negative integers and let 
$N:=k(\ell+1)+\sum_{i=1}^{\ell+1} n_i$. 
Define a multi-braid arrangement $(\A,m)=\A_\ell(n_1,\ldots,n_{\ell+1})[G]$ by $m(H_{ij}):=
2k+n_i+n_j+m_G(H_{ij})$. Assume that one of the following three conditions is satisfied:.
\begin{itemize}
\item[(1)]
$k>0$, 
\item[(2)]
$E_G^-=\emptyset$, or 
\item[(3)]
$k=0,\ m(H_{ij}) >0$ for all $1 \le i < j \le \ell+1$, and 
for all the triples $\{s,i,j\}$ with $\{i,j\} \in E_G^-$ and  
$m_G(H_{si})<m_G(H_{sj})$, it holds that $n_i>0$.

%$E_G^-=\emptyset$, or 
%\item[(3)]
%$E_G^+=\emptyset,\ m(H_{ij})>0$ for any $i \neq j$, and for any distinct triple $\{i,j,k\}$ with 
%$\{i,j\},\{i,k\} \in E_G^-$, it holds that $n_i>0$, or $\{j,k\} \in E_G^-$.
\end{itemize}
%Then $\A \in \mathcal{DF}$.
%
%derivations 
%$\theta_E, \varphi,\theta_2,\ldots,\theta_{\ell-1}$ for $D(\A)$ 
%such that $\deg \varphi=d,\ 
%\deg \theta_i=d_i\ (i=2,\lots,\ell-1)$ and $\alpha_i:=\alpha_{H_i}$ divides $\theta_i$ for $i=2,\ldots,%\ell-1$. 
Then $(\A,m)$ is free if and 
only if $G$ is signed-eliminable. In particular, when $(\A,m)$ is free, 
$\exp(\A,m)=(0,N+\widetilde{\deg}_2,\ldots,N+\widetilde{\deg}_{\ell+1})$.
\medskip

Note that the condition (1) and (2) are the same as the original conditions (i) and (ii) respectively, and 
the original proofs also work well. 
The new condition (3) corrects the original condition  (iii), and 
also contains a new multiplicity which was not in the original paper. 
By the new condition (3) in Theorem 0.3, 
the case in Example does not occur. 

To prove Theorem 0.3, let us recall the following two lemmas for the reader's convenience.
\medskip

\noindent
\textbf{Lemma A} (\cite{WA}).
Let $(\A,m)$ be a multiarrangement defined by 
$$
x^a y^b (x-y)^c=0.
$$
Assume that $a \ge b,c$. Let $d:=\lfloor (a+b+c)/2 \rfloor$. Then 

(1)\,\,
$\exp(\A,m)=(d,d)$ or $(d,d+1)$ if $a \le  b+c$. 

(2)\,\,
$\exp(\A,m)=(a,b+c)$ if $a \ge b+c+1$.
\medskip

Also, the important role was played by the Euler multiplicity and the Euler restriction 
in the paper. Let us recall it. Let $(\A,m)$ be a multiarrangement, and let $H \in \A$. Then 
we can define the Euler restriction $(\A^H,m^*)$ of $(\A,m)$ onto $H$ (see \cite{ATWA} for details), where 
$\A^H:=\{H \cap L \mid L \in \A \setminus \{H\}\}$. Since 
the Euler multiplicity $m^*$ depends only on $\A_X:=\{H \in \A \mid X \subset H\}$ and 
$m|_{\A_X}$ for $X \in L(\A)$ with $\codim X=2$, 
%the codimension two data of $(\A,m)$ 
and we are interested in the braid arrangement in the paper, 
it suffices to know the following computation for the Euler multiplicity $m^*$.
\medskip

\noindent
\textbf{Lemma B} (\cite{ATWA}). 
Let $(\A,m)$ be a multiarrangement in $\R^2$, $H \in \A$ and 
let $(\A^H,m^*)$ be the Euler restriction of $(\A,m)$ onto $H$. 

(1)\,\,
Assume that $\A=\{H,L\}$. Then $m^*(H \cap L)=m(L)$.

(2)\,\,
Assume that $\A=\{H,L,K\}$, $\exp(\A,m)=(d_1,d_2)$ and 
$\exp(\A,m-\delta_H)=(d_1,d_2-1)$, where $\delta_H$ is the characteristic multiplicity of $H$. 
Then $m^*(H \cap L)=d_1$.
\medskip

\noindent
\textbf{Remark}. Explicitly, for the braid arrangement case, the corresponding arrangement 
to those in Lemma B is the Weyl arrangement of the type $A_2$ defined by 
$$
(x-y)(y-z)(x-z)=0.
$$
Every plane of this arrangement contains a line defined by $x=y=z$, hence essentially 
in $\R^2$. Hence we may apply Lemma B to this arrangement.
\medskip

The error in the original Theorem 0.3 is based on the existence of the multiplicity for which Lemma 4.1 in the paper does not work.
Lemma 4.1 states that, for the multiarrangement in the paper, in the case of Lemma B (2), $m^*(H \cap L)=\min\{d_1,d_2\}$. By Lemma A, this 
holds true in the case (1) in Lemma A. So we have to analyze the case (2) in Lemma A.
\medskip

Hence from now on, \textbf{we consider the graphs $G$ satisfying 
one of the conditions (1), (2) or (3) in Theorem 0.3}. 

For that analysis, let us recall a definition in the paper. 
Let $G$ be a signed eliminable graph with a signed elimination ordering 
$(1,\ldots,\ell+1 )$ of $V_G$. Let $1 \le h \le \ell+1$ and let $\{v_1,\ldots,v_a\}$ be 
the vertices such that 
$v_j < h$ and  that $\{v_j,h\} \in E_G$ for all $j$. Then an $h$-th signed-eliminable filtration  (see Definition 2.2 in the paper) is 
an ordering $(i_1,\ldots,i_{a})$  of vertices $\{v_1,\ldots,v_a\}$ such that 
\begin{itemize}
\item[(1)]
if $\{h,i_s\},\{h,i_t\} \in E_G^\sigma$, then  $\{i_s,i_t \} \in E_G^{\sigma}\ (\sigma \in \{+,-\})$, and 
\item[(2)]
%the case $\{k,i_s\} \in E_G^\sigma,\ \{i_s,i_t\} \in E_G^{-\sigma}$ and 
%$\{k,i_t\} \not \in E_G$ does not occur 
if $\{h,i_s\} \in E_G^\sigma$ and $\{i_s,i_t\}
\in E_G^{-\sigma}$, then $\{h,i_t\} \in E_G^{-\sigma}$ and $t<s$ 
for $\sigma \in \{+,-\}$.
\end{itemize}

It is shown in the paper that every signed eliminable graph admits a signed eliminable 
filtration for all $h$. 
The reason why we consider a signed eliminable filtration is as follows. Assume that 
Lemma 4.1 in the original paper holds for $(\A,m)=\A_\ell(n_1,\ldots,n_{\ell+1})[G]$. Let 
$(1,\ldots,\ell+1)$ be a signed elimination ordering, $G_h$ be the induced subgraph 
from the vertices $\{1,\ldots,h\} \subset V_G$, and let 
$m_h:=m|_{G_h}$. Moreover, for the $h$-th signed eliminale filtration $(i_1,\ldots,
i_{a})$ and $1 \le s \le a$, let $G_h^s$ be the graph such that $V_{G_h^s}=V_G$ and 
$E_{G_h^s}=E_{G_{h-1}} \cup 
\{\{h,i_j\} \mid j=1,\ldots,s\}$. Define a multiplicity $m_h^s$ on $\A_\ell$ by 
$m_h^s(H_{ij}):=2k+n_i+n_j+m_{G_h^s}(H_{ij})$. 
%Then 
%it was shown in the proof of Theorem 0.3 in the original paper that 
%the Euler restriction of $(\A,m_{G_{h}^s})$ onto $H_{h,i_s}$ is
%$(\A_{\ell-1},m_h|_{\{1,\ldots,\hat{h},\ldots,\ell+1\}})$. 
%Lemma 4.1 in the original paper, and hence the addition-deletion theorem works as in the 
%paper. 
%The error in the original Theorem 0.3 is based on the case when Lemma 4.1 in the original 
%paper does not work. We have to determine in which case it occurs in 
Now we can determine when Lemma 4.1 does not work in the following lemma. 
\medskip

\noindent
\textbf{Lemma C}. 
%Let $(1,2,\ldots,\ell+1)$ be a signed elimination ordering, and let 
Under the notation above, 
assume that $G$ is signed eliminable with a signed elimination 
ordering $(1,\ldots,\ell+1)$. 
Let $(i_1,\ldots,i_a)$ be an $h$-th signed-eliminable filtration of $h,\ 1 \le h \le \ell+1$.
% i.e., 
%the induced subgraph $\{i_1,\ldots, i_k,\ell\}$ is also signed eliminable for $k=1,\ldots,s$. 
%Assume that $i,j \in \{i_1,\ldots,i_s\},\ i <j$ satisfy that 
%the Euler restriction with respect to this signed-eliminable filtration does not satisfy 
For $1 \le s \le a$ and the Euler multiplicity $m^*$ of $(\A_\ell,m_h^s)$ onto $H_{h,i_s}$, 
assume that 
$m^*$ does not satisfy 
Lemma 4.1 in the paper for $H_{h,i_s} \cap H_{h,i}$ for some $i$. 
%when restricting onto $H_{kj}$, 
Then it holds that, 
%$k=n_{i_s}=0, \{h,i_s\}, \{i_s,i\}$ and $\{h,i\} \in E_G^-$. Moreover, if $i=i_b$ in the 
%notation above, then $
either 
\begin{itemize}
\item[(a)]
$k=n_{i}=0, \{i_s,i\},\ \{h,i_s\},\{h,i\} \in E_G^-$, and 
$b<s$ for $i=i_b$, or 
\item[(b)]
$k=n_{i_s}=0, \{i_s,i\},\ \{h,i_s\},\{h,i\} \in E_G^-$, and 
$b>s$ for $i=i_b$ 
\end{itemize}
in the notation above.
\medskip

\noindent
\textbf{Proof}. For the simplicity, let $j:=i_s$. 
By Lemma A, it suffices to check the following three cases:

\textbf{Case 1}. The case when $m(H_{hj})$ is large.
In the terminology of Lemma A (2), $m(H_{hj}) \ge m(H_{hi})+m(H_{ij})+1$. However, if 
$m(H_{hj}) = m(H_{hi})+m(H_{ij})+1$, then Lemmas A and B assert that 
Lemma 4.1 in the paper holds. So we have to analyze the case 
 $m(H_{hj}) > m(H_{hi})+m(H_{ij})+1$, 
i.e., 
$$
2k+n_h+n_j+\epsilon_{hj} >4k+2n_i+n_h+n_j+\epsilon_{hi}+\epsilon_{ij}+1.
$$
Here $\epsilon_{ij}:=m(H_{ij}) \in \{-1,0,1\}$. Since we apply the addition-deletion theorems to 
$H_{hj}$, it holds that $\epsilon_{hj} \in \{0,1\}$. This inequality is  
equivalent to 
$$
\epsilon_{hj}>1+2k+2n_i+\epsilon_{hi}+\epsilon_{ij}.
$$
It is easy to check that this cannot occur unless $k=n_i=0$. In that case, this holds only when 
$$
(\epsilon_{hj},\epsilon_{hi},\epsilon_{ij})=(1,-1,-1),\ (1,0,-1),\ (0,-1,-1)\ \mbox{or}\ (1,-1,0).
$$
The case $(0,-1,-1)$ corresponds to the case (a) in Lemma C, and the other three cases do not occur by the condition in Theorem 0.3 (3) 
and the complete signed-eliminable filtration.

\textbf{Case 2}. 
The case when $m(H_{ij})$ is large, i.e., 
$$
2k+n_i+n_j+\epsilon_{ij} >4k+2n_h+n_i+n_j+\epsilon_{hi}+\epsilon_{hj}.
$$
Note that $\epsilon_{hj} \ge 0$ by the same reason as above. 
This is 
equivalent to 
$$
\epsilon_{ij}>2k+2n_h+\epsilon_{hi}+\epsilon_{hj}.
$$
It is easy to check that this cannot occur unless $k=n_h=0$ by 
the signed eliminability and the fact that $\epsilon_{hj} \ge 0$. 
In that case, again the sigied eliminability shows that this holds only when 
$$
(\epsilon_{hj},\epsilon_{hi},\epsilon_{ij})=(0,-1,0),\ (0,0,1),\ (0,-1,1)\ \mbox{or}\ (1,-1,1).
$$
These four cases do not occur by the condition in Theorem 0.3 (3)  
and the signed-eliminable filtration.

\textbf{Case 3}. The case when $m(H_{hi})$ is large, i.e., 
$$
2k+n_h+n_i+\epsilon_{hi} >4k+2n_j+n_h+n_i+\epsilon_{hj}+\epsilon_{ij}.
$$
Again $\epsilon_{hj} \ge 0$. 
This is 
equivalent to 
$$
\epsilon_{hi}>2k+2n_j+\epsilon_{hj}+\epsilon_{ij}.
$$
It is easy to check that this cannot occur unless $k=n_j=0$ by 
the signed eliminability and the fact that $\epsilon_{hj} \ge 0$. 
In that case, this holds only when 
$$
(\epsilon_{hj},\epsilon_{hi},\epsilon_{ij})=(0,0,-1),\ (0,1,0),\ (1,1,-1)\ \mbox{or}\ (0,1,-1).
$$
The case $(0,0,-1)$ corresponds to the case (b) in Lemma C, and the other three cases do not occur by the condition in Theorem 0.3 (3) 
and the signed-eliminable filtration. \owari
\medskip

Note that, by Lemma C, the error does not occur when $k>0$ or $E_G=E_G^+$, which 
corresponds to the original conditions (i) and (ii) in Theorem 0.3 in the paper.
\medskip

\noindent
\textbf{Proof of Theorem 0.3}. 
First let us show the ``if'' part. 
The original proof is correct if Lemma 4.1 works well. In other words, the original condition (iii) in Theorem 0.3 
is not enough for Lemma 4.1 to work well. Now the new condition (3) in Theorem 0.3 makes Lemma 4.1 work well in the following reason. 

%In the following assume that $k=0$. 
%If $E_G^-=\emptyset$, 
%If the case in Lemma C does not occur, then the original proof of Theorem 0.3 (ii) works well. 
%If $E_G^+ \neq \emptyset$ and $E_G^- \neq \emptyset$, then the added condition makes Lemma A (2) works %well. 
%Hence the same proof also 
%works in this case. 

By Lemma C, the original proof works well if  %$E_G^+=\emptyset$, and 
$G$ does not have a signed eliminable filtration containing 
the case (a) or (b) in Lemma C. Assume that either (a) or (b) occurs in $G$. By Lemma C and the condition 
in Theorem 0.3 (3), in this case, 
$k=0$. Hence in the following we always assume the conditon (2) in Theorem 0.3. 
%Then we need a different proof. 
We may assume that 
$(1,2,\ldots,\ell+1)$ is a signed elimination ordering for $G$. By Lemma C, 
this means that, if 
$\{\ell+1,i\}$ and $\{i,j\}$ are negative edges for $\ell+1>i,j$, and $n_i=0$, then $\{\ell+1,j\}$ is also a negative edge. 
By Case 2 in the proof of Lemma C, we may assume that $n_{\ell+1} \neq 0$, thus $n_i=0$ for some $i<\ell+1$ with 
$\{\ell+1,i\} \in E_G^-$.
We may 
easily check the statement is true for $\ell \le 3$. We use the induction on $\ell$. Assume that 
$\A_\ell(n_1,\ldots,n_{\ell+1})[G|_{\{1,\ldots,\ell\}}]$ is free with exponents 
$(0,N+\deg_1,\ldots,N+\deg_\ell,N)$.
Let us increase/decrease multiplicities on the hyperplanes $H_{j,\ell+1}$ for $j=1,\ldots,\ell$ following 
the signed eliminable filtration. 
%corresponding to edges between 
%$\ell+1$ and $1,2,\ldots,\ell$. 
%Hence for 
Let $G':=G|_{\{1,2,\ldots,\ell\}}$.
%, the multiarrangement $\A_{G'}$ is free with 
%exponents $(0,N+\deg_1,\ldots,N+\deg_\ell,N)$. 

%First, note that, 
%when $\{j,\ell+1\} \in E_G^-$ for some 
%$j \in \{1,\ldots,\ell\}$, 
%$n_{\ell+1}=n_j=0$ cannot occur. If so, then 
%$m(H_{j,\ell+1})=-1$, which contradicts the assumption on $m$.

%Hence first assume that $n_{\ell+1}>0$. 
Let $i_1<\cdots<i_a<\ell+1$ be the set of all vertices connected with $\ell+1$ by negative edges. 
%If $n_{i_j}>0$ for all $j$, then there are no problem. 
Assume that at least one of $n_{i_1},\ldots,n_{i_a}$ is zero. Then 
in fact only one of them, say $n_{i_j}$, is zero. Assume that $n_{i_j}$ and $n_{i_s}$ are both zero. Since 
$(1,\ldots,\ell+1)$ is a signed elimination ordering, $i_j$ and $i_s$ are also connected by a negative edge. Hence $m(H_{i_j,i_s})=
-1$, which is a contradiction. 
Now we show that there is an $(\ell+1)$-th signed-eliminable filtration such that ${i_j}$ is the last vertex in the filtration, i.e., 
when we apply the addition-deletion theorem for multiarrangements to 
$H_{t,\ell+1}$ for $\{t,\ell+1\} \in E_G$ following a signed-eliminable filatration, 
we may 
decrease $m(H_{\ell+1, i_j})$ in the final step of this filtration. Assume that $i_j$ is not the last vertex in the filtration, i.e., 
the filtration is of the form $(A,i_j,B)$ for ordered sequences $A$ and $B$ of vertices connected with $\ell+1$.
Then consider the new filtration $(A,B,i_j)$. We show that this is also an $(\ell+1)$-th
signed-eliminable filtration. Assume not. Then 
there is $h,\ 1 \le h \le \ell$, such that $\{\ell+1, i_j\} , \{i_j, h\}\in E_G^-$ and that $\{\ell+1,h\} \in E_G^+$. 
%for some $h<\ell+1$ with 
However, this 
says that $n_{i_j} \neq 0$ by the assumption in Theorem 0.3 (3), 
which 
is a 
contradiction.

So we may use the same argument as in the original paper to 
increase/decrease multiplicities $m(H_{t,\ell+1})$ for $i_j \neq t <\ell+1$ with $\{t,\ell+1\} \in E_G$ to obtain a free arrangement with 
exponents $(0,N-\deg_1,\ldots,N-\deg_\ell,N-\deg_{\ell-1}+1)$. In the final step decrease the multiplicity $m(H_{i_j, \ell+1})$ by one. Since all the other vertices connected to $\ell+1$ 
are increased/decreased, we may show that also here the same argument as in the original paper works by Lemmas A and B, which completes the proof. 

Also, note that the Euler restriction along an $h$-th signed eliminable filtration does not change 
the edges $\{i,j\}$ and its sign $\pm$ for $i,j<h$ as we saw above. Hence if $G$ satisfies the condition (3) in Theorem 0.3, 
then so is the restricted graph, which makes the induction work. 

%Note that the proof above confirms that if one of the conditions (1), (2) or (3) is satisfied, then 
%Lemma 4.1 in the original paper holds. 

Next we show the ``only if'' part. 
%The original proof works 
%when the conditions (1) or (2) is satisfied. Consider 
%the case when the condition (3) holds. 
The proof of ``only if'' part consists of two arguments in the original paper.
The first one is in page 10 in the original paper, the inequalities $B_i$. 
In this case, the proof is the same 
as the original one. What we have to pay attention is whether $\exp(\A_X,m_X)$ 
is of the form in Lemma A (1) for all codimension two flat $X$.
This is confirmed by the conditions (1), (2) and (3) in Theorem 0.3.
%, and 
%the Euler multiplicity satisfies Lemma 4.1 in the original paper. 
The second one is to investigate when $G$ contains a cycle of length at least four, 
mountains or hills defined in the original paper. To apply these arguments, we need Lemma 
4.1 to apply the addition-deletion theorems. Now this works well by Lemma C and the condition (3) in Theorem 0.3.
%if 
%$\exp(\A_X)=(e_1,e_2,0,\ldots,0)$ with $|e_2-e_1|\le 1$ for all $X \in L(\A)$ with $|\A_X|=3$. By Lemmas A and %B, we can check that this is true when the condition 
%(1) or (2) in Theorem 0.3 is satisfied. 
%\owari
\owari
%
%Next assume that $n_{\ell+1}=0$. But this case was already removed in Case 2. \owari
\medskip

\noindent
\textbf{Remark}. The same statement holds true even if $m(H_{ij})=0$ could occur.
\medskip

Except for the condition (iii) in Theorem 0.3, all the results and proofs are correct in this paper. 
\medskip

\noindent
\textbf{Acknowledgements}. The authors are really grateful to Michael Dipasquale for his 
pointing out the mistake of the original Theorem 0.3 with a counter example, and his careful reading 
of this draft with useful comments.

\end{document}